\documentclass[12pt]{article}

\newcommand{\ben}{\begin{enumerate}}
\newcommand{\een}{\end{enumerate}}
\newcommand{\be}{\begin{equation}}
\newcommand{\ee}{\end{equation}}
\newcommand{\bse}{\begin{subequation}}
\newcommand{\ese}{\end{subequation}}
\newcommand{\bea}{\begin{eqnarray}}
\newcommand{\eea}{\end{eqnarray}}
\newcommand{\bc}{\begin{center}}
\newcommand{\ec}{\end{center}}

\def\b#1{{\bf #1}}
\def\c#1{{\cal #1}}
\def\1{{\bf 1}}
\def\ket#1{\vert #1 \rangle}

\author{
        {\bf  Gaetano Fiore}\\
        {\it Dip. di Matematica e Applicazioni, Fac.  di Ingegneria,}\\ 
        {\it Universit\`a di Napoli, V. Claudio 21, 80125 Napoli}\\
        {\it and}\\
        {\it I.N.F.N., Sezione di Napoli,
        Mostra d'Oltremare, Pad. 19, 80125 Napoli}\\
        \\
        {\bf  John Madore}\\
        {\it Laboratoire de Physique Th\'eorique et Hautes Energies,}\\
        {\it Universit\'e de Paris-Sud, B\^atiment 211, F-91405 Orsay}\\
        }
\title{Geometrical Issues for the 3-dim Quantum Euclidean Space
\footnote{Talk given by the first author at the
``VI Wigner Symposium'', Istanbul, August 1999}}

\begin{document}
\maketitle

\begin{abstract}
We briefly describe our application of a version of
noncommutative differential geometry to the 3-dim quantum
space covariant under the quantum group of rotations $SO_q(3)$ and
sketch how this might be used to determine the correct
physical interpretation of the geometrical observables.
\end{abstract}

\vskip1cm
\bc
Preprint 99-53, Dip. Matematica e Applicazioni, Universit\`a di Napoli
\ec

\section{Introduction and preliminaries}

It is a rather old idea that the
micro-structure of space-time at the Planck level might be better
described using a noncommutative geometry. 
An
interesting problem is that of finding an appropriate noncommutative
version of Minkowski space (see e.g. ~\cite{Wess,Castellani}). 
Here we consider a
noncommutative generalization~\cite{DimMad96} of the Cartan
`moving-frame' formalism which provides in certain special cases an
interesting bridge between the `Dirac-operator' formalism of
Connes~\cite{Con94} and the quantum-group formalism of Woronowicz
\cite{Wor89}. We apply it to the quantum Euclidean space 
$\b{R}_q^3$ \cite{FRT}, namely the quantum space covariant under 
the quantum group $SO_q(3)$. We first outline
the main results of a previous article~\cite{FioMad98'} and
then sketch how these results can be used to extract informations on the
geometrical structure of $\b{R}_q^3$. For the generalizations of these
results to the quantum Euclidean spaces of higher dimensions see
Ref. \cite{CerFioMad2000} and the article entitled ``Geometrical 
Techniques on the $N$-dimensional Quantum Euclidean Spaces'' in 
the present proceedings. 
In Ref. \cite{FioMad98'}
we introduced a metric and an `almost' metric-compatible
linear connection on the quantum Euclidean 
space, equipped with its (two) standard $SO_q(3)$-covariant differential
calculi; correspondingly, the `frame' or dreibein has been also found. 
Modulo a conformal factor, which might however be
reabsorbed into a formulation of the metric
compatibility more suitable for the present case, the curvature 
turns out to be zero, suggesting that the quantum space is flat as in the 
commutative limit. This is of course welcome if the postulated
noncommutative algebra and differential calculus are really to
describe flat Euclidean space $\b{R}^3$
in the commutative limit. In a separate
paper we shall show that in the same limit the traditional quantum space 
coordinates go to suitable general (non-cartesian) coordinates. 
This will allow to cure some unpleasent features \cite{fiolat} of a naive 
physical interpretation of the representation theory of 
the algebra of function on $\b{R}_q^3$.
Using the properties of the volume form on $\b{R}_q^3$, here we just
give some arguments why this phenomenon must occur. 

\medskip

The preliminaries contained in this section are 
a variation of noncommutative geometry which has been 
proposed~\cite{DubMadMasMou95, DimMad96} as a noncommutative version
of the Cartan `moving-frame' calculus.
The starting point is a noncommutative algebra $\c{A}$ with
unit which has as commutative limit 
the algebra of functions on some manifold $\c{M}$ and over
$\c{A}$ a differential calculus~\cite{Con94} 
$\{d,\Omega^*(\c{A})\}$ which 
we shall choose so that it
has as corresponding limit the ordinary de~Rham differential
calculus. It is determined 
by the left and right module structure of the
$\c{A}$-module of 1-forms $\Omega^1(\c{A})$.  
By definition a {\it metric} is a nondegenerate 
$\c{A}$-bilinear map
\be
g:\Omega^1(\c{A})\otimes_{\c{A}} \Omega^1(\c{A})\rightarrow
\c{A}.
\ee
$\c{A}$-bilinearity means
\be
g(f\xi\otimes\eta h)=f g(\xi\otimes\eta) h,
\ee
for any $f,h\in\c{A}$ and $\xi, \eta\in\Omega^1(\c{A})$. 
It implies
that $g$ is completely determined by the $\c{A}$-valued 
matrix elements
\be
g^{ab}:=g(\theta^a\otimes\theta^b),        \label{defgab}
\ee
for any choice of a basis of 1-forms $\{\theta^a\}$.
In the commutative limit $\c{A}$-bilinearity is 
equivalent to the very important requirement of locality  
of $g$ in both arguments at each point $x\in\c{M}$:
\be
[g(f\xi\otimes\eta h)](x)=f(x)\: [g(\xi\otimes\eta)](x)\: h(x).
\ee

A {\it linear connection} is a map (see \cite{Kos60})
\be
D:\Omega^1(\c{A}) \rightarrow\Omega^1(\c{A})\otimes_{\c{A}} 
\Omega^1(\c{A})
\ee
together with a ``generalized flip'' $\sigma$, i.e. a 
$\c{A}$-bilinear map
\be
\sigma:\Omega^1(\c{A})\otimes_{\c{A}} \Omega^1(\c{A}) 
\rightarrow\Omega^1(\c{A})
\otimes_{\c{A}} \Omega^1(\c{A})
\ee
going to the ordinary flip in the commutative limit and
such that $D$ satisfies the left and right Leibniz rules
\bea
D (f \xi)  & =&  df \otimes \xi + f D\xi                    
\label{2.2.2}\\
D(\xi f) &= &\sigma (\xi \otimes df) + (D\xi) f .          
\label{second}
\eea
Because of bilinearity,
given a basis of 1-forms $\{\xi^i\}$ $\sigma$ is completely 
determined once the
$\c{A}$-valued matrix elements $S^{ij}{}_{hk}$ defined by
\be
\sigma_0 (\xi^i\otimes \xi^j)  =: S^{ij}{}_{hk}\,
\xi^h\otimes \xi^k
   \label{assign}
\ee
are assigned. 

Let $\pi$ be the projection
\be
\pi:\Omega^1(\c{A})\otimes_{\c{A}} \Omega^1(\c{A}) 
\rightarrow\Omega^2(\c{A}).
\ee
The {\it torsion} is the map $\Theta= d-\pi\circ D$. 
In order 
that the torsion be bilinear we shall require
\be
\pi\circ (\sigma+ \1)=0.
\ee
One can naturally extend $D$ to higher tensor powers, e.g.
\be
D_2(\xi\otimes\eta)=
D\xi\otimes\eta+\sigma_{12}(\xi\otimes D\eta),
\ee
where we have introduced the tensor notation 
$\sigma_{12}=\sigma\otimes \1$.
The metric-compatibility condition for $g,D$ 
reads $g_{23}\circ D_2=d\circ g$.

The {\it curvature} 
$\mbox{Curv}:\Omega^1(\c{A})\rightarrow\Omega^2(\c{A})
\otimes_{\c{A}} \Omega^1(\c{A})$ is defined by
\be
\mbox{Curv} = \pi_{12}\circ D_2\circ D.         \label{curv}
\ee
It is always left $\c{A}$-linear, and right $\c{A}$-linear 
only in certain
models; in general, right linearity is guaranteed only in 
the commutative
limit. Therefore in this limit the curvature is local,  an 
essential physical requirement for a reasonable definition 
of a curvature.

If $\c{A},\Omega^1(\c{A})$ are $*$-algebras and $d$ is real,
$(df)^*=df^*$,  $D$ is said to be real  \cite{FioMad98} if
\be
D\xi^*=(D\xi)^*
\ee
where the involution
on $ \Omega^1(\c{A})\otimes_{\c{A}} \Omega^1(\c{A})$ is defined by 
\be
(\xi\otimes\eta)^*=\sigma(\eta\otimes\xi^*),
\label{star1}
\ee
with a $\sigma$ such that the square of $*$ gives the identity. 
Note that this expression has the correct classical limit. So
real structures on the
tensor product are in one-to-one correspondence with right 
Leibniz rules.
$D_2$ is real if $D$ is and $\sigma$ in addition fulfills 
the braid equation
\be
\sigma_{12}\sigma_{23}\sigma_{12}=\sigma_{23}\sigma_{12}
\sigma_{23}.                                                  
\label{genbraid}
\ee
The curvature is real if $D$ is real and
(\ref{star1}), (\ref{genbraid}) are satisfied.

Now assume that there exists a {\it frame}, i.e.
a special basis $\theta^a\in\Omega^1(\c{A})$,  
$1\leq a \leq n$, 
such that
\be
[\theta^a,\c{A}]=0 
\ee
and any $\xi\in\Omega^1(\c{A})$ can be uniquely written 
in the form
$\xi_a\theta^a$, with $\xi_a\in\c{A}$. This is possible only
if the limit manifold $\c{M}$ is parallelizable. It has 
the advantage 
that for any $f\in\c{A}$ the computation of
commutator $[\xi,f]$ is reduced to the computation of
the commutators $[\xi_a,f]$ in $\c{A}$.
Assume also that there exist $n$ inner derivations $e_a$, 
\be
e_a f:=
[\lambda_a, f]
\ee
($\lambda_a\in\c{A}$), dual to $\theta^a$: 
$\theta^a(e_b)=\delta^a_b$. Then
\be
\theta := - \lambda_a \theta^a                                 
\label{dirac}
\ee 
is the `Dirac operator'~\cite{Con94} for d:
\be
df = -[\theta,f].                                              
\label{extra}
\ee
$\theta^a$ is a very convenient basis to work with. 
For instance, from 
$\c{A}$-bilinearity it immediately follows that  
the corresponding
elements (\ref{defgab}) lie in the center 
$\c{Z}(\c{A})$ of $\c{A}$,
by the sequence of identities
\be
f g^{ab} = g(f \theta^a \otimes \theta^b) 
= g(\theta^a \otimes \theta^b f) = g^{ab} f.              
\label{2.2.22}
\ee 
We shall be interested in the
case that $\c{Z}(\c{A})={\bf C}$. In the commutative limit
the condition $g^{ab}\in{\bf C}$ characterizes the vielbein or  
`moving frame' of E. Cartan, which is determined
up to a linear transformation; 
if this condition is fulfilled for any 
value of the deformation parameter the $\theta^a$ remain
uniquely determined up to a linear transformation
and are particularly convenient objects to be used to guess 
a physically sensible formulation of noncommutative-geometric 
notions. 

{From} the above considerations we deduce that the metric is 
fixed by
the form of the frame, and since the latter is fixed by the 
structure
of the differential calculus up to
a $GL(n)$ transformation, it will also be. Hence the differential 
calculus will yield a metric in the commutative limit as a
{\it shadow of noncommutativity}. As a consequence, it seems 
that in this framework the metric cannot be considered as a dynamical 
variable, since it is completely determined 
by the differential calculus. This is completely different
from what occurs in the commutative case, where the
differential calculus and the metric on a smooth
manifold are two independent structures. 
Something basically similar was proposed
many years ago by Wheeler~\cite{MisThoWhee73} when he
suggested that the `graviton' be considered as the `phonon' of a
fundamental space-time lattice. The fact that Wheeler was 
considering an
ordinary lattice instead of a `quantum lattice' is important 
of course but
not essential.

One might be tempted to use this fact 
as an argument against noncommutative geometry.
However one should note that in NCG the freedom lost in choosing
the metric is recovered as a much larger freedom in choosing
the differential calculus (without necessarily changing the 
commutative limit of the latter). In other words, the 
`degrees of freedom' of the metric will be now encoded in the
structure of the differential calculus.

\section{Application of the formalism to the quantum Euclidean space}

Take  `the algebra of functions on the quantum Euclidean
space $\b{R}_q^3$' 
\cite{FRT} as $\c{A}$ and over it one of the two $SO_q(3)$-covariant
differential calculi \cite{cawa}. The treatment
of the other calculus can be done in a completely
parallel way, see ref. \cite{FioMad98'}.
Here we are interested in the case of
a real positive $q$. We ask if they fit in the
previous scheme. We shall denote by $\Vert \hat R^{ij}_{hk}\Vert$ 
the braid matrix
of $SO_q(3)$, by $g_{ij}=g^{ij}$ the $SO_q(3)$-covariant
metric; here and below all indices will take the
values $-,0,+$. In the commutative limit $q\rightarrow 1$
$g_{ij}\rightarrow \delta^{i,-j}$.
The projector decomposition of $\hat R$ is
\be
\hat R = q\c{P}_s - q^{-1}\c{P}_a + q^{-2}\c{P}_t ;            
\label{decompo}
\ee
$\c{P}_s$, $\c{P}_a$, $\c{P}_t$ are $SO_q(3)$-covariant
$q$-deformations of respectively the symmetric trace-free,
antisymmetric and trace projectors. The trace projector is
1-dimensional and is related to $g_{ij}$ by
\be
\c{P}_t{}_{kl}^{ij} \propto g^{ij}g_{kl},         \label{bibi}
\ee
$\c{A}$ is generated by $x^-,x^0,x^+$ fulfilling $\c{P}_a xx=0$,
or more explicitly
\bea
&&x^- x^0 = q\, x^0 x^-,\nonumber\\
&&x^+ x^0 = q^{-1} x^0 x^+, \label{xxcr3}\\
&&[x^+, x^-] = h (x^0)^2.\nonumber
\eea
where we define $h = \sqrt q - 1/\sqrt q$. For real positive $q$
the real structure on $\c{A}$
is defined by $(x^i)^*=x^jg_{ji}$, or more explicitly
\be
(x^-)^* = \sqrt q x^+, \qquad (x^0)^* = x^0, \qquad
(x^+)^* = 1/\sqrt{q} x^-.
\ee
$\c{Z}(\c{A})$ is generated by the $SO_q(3)$-covariant real element
\be
r^2 := g_{ij} x^i x^j = \sqrt{q}x^+x^-+ (x^0)^2+1/\sqrt{q}x^-x^+.
\label{sql}
\ee
Let $\xi^i = dx^i$. One $SO_q(3)$-covariant calculus, which we shall
denote by $\{d,\Omega^*(\c{A})\}$, is determined by the commutation
relations
\be
x^i \xi^j = q\,\hat R^{ij}_{kl} \xi^k x^l.                     
\label{xxicr}
\ee
Unfortunately for real positive $q$
neither calculus has a real exterior derivative, 
and up to now no way was known to make it closed under involution
\cite{olezu}; rather, each exterior algebra is mapped into
the other under the natural involution.
The `Dirac operator' (\ref{extra}) corresponding to $d$ is the
$SO_q(3)$-invariant element \cite{zumino}
$\theta := (q-1)^{-1} q^2 r^{-2} x^i\xi^j g_{ij}$;
note that $\theta$ is singular in the commutative limit.
As a contribution to the understanding of the structure of the quantum
Euclidean spaces we have noticed~\cite{FioMad98'} the following results

\begin{enumerate}
\item There exist two torsion-free, `almost' metric-compatible
linear connections, given by the formula
\be
D_{(0)} \xi = -\theta \otimes \xi + \sigma_{0}
(\xi \otimes \theta)                            \label{conn}
\ee 
The two corresponding generalized flips $\sigma_{0}$ 
are the ones with matrix (\ref{assign}) given
respectively by $S=q\hat R,(q\hat R)^{-1}$.
$D_{(0)}$ `almost' metric-compatible means compatible 
up to a conformal factor with the metric which we
shall give in the next item;
a strict compatibility does not seem possible.
Both $\sigma_{0}$ fulfill the braid equation (\ref{genbraid})
and both $D_{(0)}$ are $SO_q(3)$-invariant.
\item If we extend $\c{A}$ by adding the `dilatation' generator 
$\Lambda$  
\be
x^i\Lambda= q\Lambda x^i  
\ee
together with its inverse $\Lambda^{-1}$ (we shall normalize them 
so that $\Lambda^*=\Lambda^{-1}$) and set $d\Lambda=0$, then 
up to normalization there exists a unique metric $g_0$, 
\be
g_0(\xi^i\otimes \xi^j)  = g^{ij}\,r^2\Lambda^2
\ee
($g_{ij}$ is the $SO_q(3)$-covariant metric matrix),
which is compatible with the two $D_{(0)}$ up to the conformal
factors $q^2,q^{-2}$,
\be
S^{ij}{}_{hk}g^{kl}S^{mn}{}_{jl}=q^{\pm 2}g^{im}\delta^n_h,
\ee
respectively in the cases $S=(q\hat R)^{\pm 1}$. A strict 
compatibility
would have required no $q^{\pm 2}$ at the rhs.
\item Curv=0 for both $D_{(0)}$.
\item If we further extend $\c{A}$ by adding also the generators
$r$ [the square root of (\ref{sql})], its inverse $r^{-1}$ and
the inverse $(x^0)^{-1}$ of $x^0$, then there exist a frame
$\theta^a$ , $a=-,0,+$,
and a dual basis $e_a$ of inner derivations given by 
\be
\theta^a := \Lambda^{-1}\, \theta^a_i \xi^i            \label{ok}
\ee
with
\be
\Vert\theta^a_i\Vert:= \left\Vert
\begin{array}{ccc}
(x^0)^{-1}                            &                  & \cr
\sqrt{q}(q+1) (rx^0)^{-1}x^+          & r^{-1}           &\cr
-\sqrt{q}q(q+1) (r^2x^0)^{-1}(x^+)^2  & - (q+1)r^{-2}x^+ 
& r^{-2}x^0 \cr
\end{array}
\right\Vert
\ee
\be
\begin{array}{l}
\lambda_- = + h^{-1} q \Lambda(x^0)^{-1} x^+,         \\   
\lambda_0 = - h^{-1} \sqrt q \Lambda (x^0)^{-1} r,\\
\lambda_+ = - h^{-1} \Lambda (x^0)^{-1} x^-.
\end{array}                                                     
\label{lambda'}
\ee
$e_a x^i= q\Lambda \, e^i_a$, where $\Vert e^i_a\Vert$ is 
(left and right) inverse of the $\c{A}$-valued matrix 
$\Vert\theta^a_i\Vert$.
Its elements fulfill the `$RTT$-relations' \cite{FRT}
\be
\hat R_{kl}^{ij}\, e^k_a e^l_b
= e^i_c e^j_d \,\hat R_{ab}^{cd}                            
\label{RTT}
\ee
as well as the `$gTT$-relations'
\be
g^{ab} e^i_a e^j_b = r^2 g^{ij}\qquad\qquad
g_{ij} e^i_a e^j_b = r^2 g_{ab}                             
\label{gTT}.
\ee
In a sense $r^{-1}e^i_a$ are a realization of the generators
$T^i_a$ of $SO_q(3)$.
As a consequence we find
\be
\c{P}_t{}^{ab}_{cd} \theta^c \theta^d = 0 \qquad\qquad
\c{P}_s{}^{ab}_{cd} \theta^c \theta^d = 0,               
\label{ththcr}
\ee
the same commutation relations fulfilled by the $\xi^i$'s. 
Similarly, the $\lambda_i$ and the $x^i$ satisfy the same 
commutation 
relations. Finally, up to a normalization
$g_0(\theta^a\otimes\theta^b) = g^{ab}$ .
\item $\Omega^*(\c{A})$ is closed under the involution defined by
\be
(x^i)^{\star}= x^jg_{ji} \qquad\qquad
(\theta^a)^{\star} =\theta^b g_{ba} 
\ee
(the latter acts nonlinearly on the $\xi^i$'s: 
$(\xi^i)^{\star}=\Lambda^{-2} \xi^j\,c_{ji}$,
with non-constant $c_{ji}\in\c{A}$).
\end{enumerate}

The reality structure of these differential calculi is an old, 
well-known problem (see \cite{olezu}). The solution
proposed in item 5 is not fully satisfactory, at least naively.
For instance, it does not yield real $d,D$; only the curvature 
is real, for the simple reason
that it vanishes. The involution cannot be consistently extended
to $\Omega^*(\c{A})\otimes\Omega^*(\c{A})$ according
to (\ref{star1}). Finally, apparently it has not
the correct classical limit. Actually, the latter point can be solved 
by a more careful analysis \cite{FioMad98'} leading to the
identification of $x^i$ with some suitable general coordinates, as
e.g. the ones reported at formulae (\ref{special}) below.

A more careful analysis is needed at this point, but is out
of the scope of the present report (for more details
see Ref. \cite{FioMad98'}). It involves the investigation
of the properties of the $*$-representations of $\Omega^*(\c{A})$
and seems to suggest a
more sophisticated version of the proposal in item 5, in which
the opposite properties of the two differential calculi
cancel with each other. The problems mentioned above and the 
fact that the
linear connections $D_{(0)}$ are metric-compatible up
to conformal factors (or, in other words, are only conformally flat)
may be related, in the sense that a satisfactory formulation
of the reality properties could eventually yield also a new and
satisfactory 
formulation of metric-compatibility which can be strictly fulfilled.
A careful analysis of the commutative limit is also needed 
in order to propose a reasonable correspondence principle 
between the
`new' theory and classical differential geometry. 

\section{Representation theory and geometry}

The set $(x^i, e_a)$ generates the phase space algebra 

$\c{D}_h$ of 
observables of
a `point particle on $\b{R}_q^3$'; we shall assume that
$x^i$ generate the subalgebra $\c{A}_p$ of position (i.e. 
configuration space)
observables. In order to understand the geometrical structure 
of configuration space one should first consider
the irreducible 
$*$-representations of $\c{D}_h$ on Hilbert spaces $H$, and then
attach to the configuration space observables $x^i$ the 
appropriate
physical meaning, with the help of the geometric tools (metric,
curvature, etc) described in the previous sections.
This will be done in detail elsewhere. Here we just give 
a flavour of
how this may drastically change our naive expectations 
about the
physical meaning of the generators $x^i$ of $\c{A}$, 
in the sense that
they should be interpreted as a noncommutative generalization
of generalized rather than cartesian coordinates on 
flat $\b{R}^3$.
This will automatically cure some unpleasent features 
\cite{fiolat}
which make the physical interpretation of $x^i$ as 
cartesian coordinates
problematic, namely that within each irreducible $*$ 
representation
the spectrum of $x^0$ has all eigenvalues of the same sign, 
that 
the value $0$ is an accumulation point of the spectrum of 
$|x^0|$
on the left, whereas the difference between two neighbouring 
points of
this spectrum diverges for $|x^0|\to\infty$.

A complete set of independent commuting observables (CSICO) 
must have three
elements (three being the dimension of the classical 
underlying manifold).
One cannot find three such elements within the 
subalgebra $\c{A}_p$,
since the latter is
not abelian.
As a CSICO we choose \cite{fiolat}
\be
r,x^0,{\bf k},
                                     \label{observables}
\ee
where ${\bf k}$ can be identified with $q^{L_0}$ and $L_0$ 
with an ordinary
angular momentum component along the direction of an axis 
$y^0$ in
ordinary Euclidean space (so the spectrum of $L_0$ is $\b{Z}$).
Assume that $q>1$,
for the sake of being specific.  In Ref.'s \cite{fiojmp,fiolat} 
it was shown
that the irreducible $*$-representations of $\c{D}_h$ 
for a zero spin 
point particle are essentially parametrized by
a sign $\eta=\pm$ and a constant $c\in[1,q)$, and that 
for each
$(\eta,c)$ there exists a corresponding orthonormal basis
$\{\ket{n_0,n,m}\}$ in which $x^0,r$ are diagonal:
\bea
r\ket{n_0,n,m} &&= c\,q^{n}\ket{n_0,n,m}, \label{rep2first}   
\\[2pt]
x^0 \ket{n_0,n,m} &&= \eta\,c\,
q^{(n - n_0 - \frac 12)}\ket{n_0,n,m}, \\[2pt]
\Lambda\ket{n_0,n,m} &&= \ket{n_0,n+1,m}. \label{rep2last}
\eea
The integers are such that $n_0$ in $\b{N}\cup \{0\}$ 
and $m$, $n$ in $\b{Z}$.
$m$ can be identified with the eigenvalue of the angular 
momentum component 
along a cartesian coordinate $y^0$. 

Here we do not write down the explicit action of the 
remaining generators
of $\c{D}_h$ on this basis. We just note that by 
applying both sides of 
the identity $e_af=[ \lambda_a,f]$ to the generic vector 
$|\psi\rangle\in H$ one finds that it is consistent to set
\be
e_a |\psi\rangle= \lambda_a |\psi\rangle,               
\label{lalla}
\ee
where at the rhs the action of the element $\lambda_a $ on $H$  
must be understood.

Assume for one moment that one could represent the exterior 
algebra 
$\Omega(\c{A})$ on the same Hilbert space $H$ (as we shall see, 
strictly speaking this is not possible). Since
$\theta^a$ commutes with $\c{A}$, from (\ref{lalla}) we see
that as an operator on $H$ $\theta^a$ commutes with the whole 
$\c{D}_h$,
\be
[\theta^a,\c{D}_h]=0.
\ee
Hence it is a {\it Casimir} of the representation:
\be
\theta^a |\psi\rangle= t^a|\psi\rangle 
\ee
for any $|\psi\rangle \in H$, where the objects $t^a$ 
characterize
$H$. Also wedge products of $\theta^a$  will commute 
with the whole 
$\c{D}_h$, hence
\bea
&&\theta^a \theta^b|\psi\rangle= t^at^b|\psi\rangle   \\
&&\theta^a \theta^b \theta^c|\psi\rangle\equiv\epsilon^{abc} 
\theta^+ \theta^0 \theta^-= \epsilon^{abc} (t^+t^0t^-)|\psi\rangle   
\eea
where $\epsilon^{abc}$ is the $q$-epsilon tensor \cite{fio94}.
Strictly speaking,
one cannot take $t^a$ in $\b{C}$ because then the $t^a$ would
commute with each other and therefore would not respect
the commutation relations (\ref{ththcr}). One can however assume
that at least $dv:= t^+t^0t^-$ is a constant. Using the
commutation relations and the values of $\epsilon^{abc}$ one 
can easily show that $dV:= \theta^+ \theta^0 \theta^-$ is real 
w.r.t. the star structure $\star$; so $dv\in\b{R}$.

What is the physical interpretation of $dv$, the (unique) 
eigenvalue
of the volume element observable
$dV$? It seems natural to consider it as the
volume of the generic elementary cell in the configuration space 

lattice. So the latter will be the same all over the configuration 
space. This is welcome because it means that the
uncertainty in the localization of a point particle will be 
essentially 
the same all over the space. 

One might argue that there can be no elementar cell in 
configuration
space, since we cannot choose a CSICO consisting of three
elements of $\b{R}_q^3$. With the choice (\ref{observables}) 
there are 
just
elementar ``annuli''  $V_{n,n_0}$, defined by 
$cq^n\le r<cq^{n+1}$, 
$q^{n_0}\le \frac{|x^0|}r<q^{n_0+1}$.
So in a sense we can just argue that the integral 
of $dV$ on such
annuli will give the volume of the latter. However 
this can
be computed by a regularized trace on the eigenvectors of $L_0$, by
\be
V_{n,n_0}=C\lim\limits_{N\to \infty}\frac 1{2N+1}\sum\limits_{m=-N}^N
\langle n,n_0,m|dV|n,n_0,m\rangle=C\lim\limits_{N\to \infty} 
dv=C \,dv
\ee
($C$ is an arbitrary normalization constant), which is 
also a constant.
We shall now use this result to show that the generators 
$x^i$ cannot
go to cartesian coordinates on $\b{R}^3$ in the commutative 
limit.
In the commutative case the shell $R\le r<R+\Delta R$, 
with $r$
defined by $r=\sqrt{y\cdot y}$ and $y^i$ cartesian coordinates,  
is a 
spherical shell in $\b{R}^3$, and therefore has a finite volume. 
On the contrary, in the present
noncommutative space the shell $cq^n\le r<q^{n+1}c$, with $r$
defined by $r=\sqrt{x\cdot x}$, will have an infinite volume, 
since
it can be divided in an infinite number
of annuli as above, with $n_0\in\b{N}\cup \{0\}$, each of which
has the same volume. Obviously 
this result remains true when we take the commutative
limit. Therefore in this
limit the $x^i$ might only go to some generalized, rather
than cartesian, coordinates. Let 
\be
x^0=f(\vec{y})     \qquad r= g(\vec{y})               \label{lulu}
\ee
be the transformation from some cartesian coordinates $y^i$ on 
$\b{R}^3$
to the commutative limit of the $x^i$. We can give a 
sense to this map also when $q\neq 1$, since we have assumed 
the commuting 
generators $x^0,r$ to be position observables for any $q$.
In Ref. \cite{FioMad98'} we have analyzed in some detail the 
constraints 
that a number of formal requirements puts on the commutative 
limit. 
One possible solution to these constraints gives for the 
functions $f,g$ of
(\ref{lulu})
\be
\begin{array}{l}
x^0= e^{\alpha y^0-\frac{\alpha^3}2}  \\
r=e^{-\alpha^3+\frac{\alpha^2}2y^+y^-+\alpha y^0} 
\end{array}
\label{special}
\ee
with $e^{\alpha^3}=q$. Thus, the surfaces $r=const$ are 
interpreted in physical 
space as paraboloids with axis $y^0$ rather than spheres with 
center in the 
origin, the
surfaces $x^0=const$ are interpreted in physical 
space as planes perpendicular to $y^0$ (exactly as before), 
the surfaces
$x^0/r=const$ are interpreted in physical 
space as cylinders with axis $y^0$ rather than as cones 
centered at the origin, 
and the lines $x^0=const$, $r=const$ are interpreted in physical 
space as circles perpendicular to and with center on the 
axis $y^0$.  
The exponential relation between $x^0$ and $y^0$ is analogous 
to the one
found \cite{CerHinMadWes98} for a 1-dimensional $q$-deformed 
model.
Due to the quantization of $x^0,r$, also
$y^0$ and $y_{\perp}=\sqrt{y^+y^-}$ are quantized. $y^0$ will take 
the values $y^0=\alpha p$, whereas $y_{\perp}=\alpha \sqrt{n_0}
+\frac 12$. So the eigenvalues of $y^0$ are equidistant and both 
positive and negative. Also, note that the step $\alpha$  
goes to zero when $q\to 1$. 
Thus, the unpleasent features mentioned
in the first paragraph of the present sections have been cured.

\end{document}